\documentclass[11pt]{article}
\usepackage{amsbsy}
\usepackage{amsmath, amsthm, amssymb}
\usepackage{makeidx}

\usepackage{graphicx}
\usepackage[colorlinks=true]{hyperref}
\begin{document}

\title{\bf Flexible Cholesky GARCH model with time dependent coefficients}

\vspace{3.5cm}

\author{{Toktam Valizadeh$^1$, Saeid Rezakhah$^1$}\vspace{.5cm}
\\\it
\small $^1$Department of Statistics, Faculty of Mathematics and Computer Science,\\
\small Amirkabir University of Technology, Tehran, Iran.
}
\date{}
\maketitle

\begin{quotation}
\noindent {\it Abstract:}
\par
\end{quotation}\par

\newcommand{\balpha}{\boldsymbol{\alpha}}
\newcommand{\bbeta}{\boldsymbol{\beta}}
\newcommand{\bgamma}{\boldsymbol{\gamma}}
\newcommand{\bdelta}{\boldsymbol{\delta}}
\newcommand{\bepsilon}{\boldsymbol{\epsilon}}
\newcommand{\bvarepsilon}{\boldsymbol{\varepsilon}}
\newcommand{\bzeta}{\boldsymbol{\zeta}}
\newcommand{\bet}{\boldsymbol{\eta}}
\newcommand{\btheta}{\boldsymbol{\theta}}
\newcommand{\biota}{\boldsymbol{\iota}}
\newcommand{\bkappa}{\boldsymbol{\kappa}}
\newcommand{\blambda}{\boldsymbol{\lambda}}
\newcommand{\bmu}{\boldsymbol{\mu}}
\newcommand{\bnu}{\boldsymbol{\nu}}
\newcommand{\bxi}{\boldsymbol{\xi}}
\newcommand{\0}{\boldsymbol{0}}
\newcommand{\1}{\boldsymbol{1}}
\newcommand{\bpi}{\boldsymbol{\pi}}
\newcommand{\bvarpi}{\boldsymbol{\varpi}}
\newcommand{\brho}{\boldsymbol{\rho}}
\newcommand{\bvarrho}{\boldsymbol{\varrho}}
\newcommand{\bsigma}{\boldsymbol{\sigma}}
\newcommand{\bvarsigma}{\boldsymbol{\varsigma}}
\newcommand{\btau}{\boldsymbol{\tau}}
\newcommand{\bupsilon}{\boldsymbol{\upsilon}}
\newcommand{\bphi}{\boldsymbol{\phi}}
\newcommand{\bvarphi}{\boldsymbol{\varphi}}
\newcommand{\bchi}{\boldsymbol{\chi}}
\newcommand{\bpsi}{\boldsymbol{\psi}}
\newcommand{\bomega}{\boldsymbol{\omega}}
\newcommand{\bGamma}{\boldsymbol{\Gamma}}
\newcommand{\bDelta}{\boldsymbol{\Delta}}
\newcommand{\bTheta}{\boldsymbol{\Theta}}
\newcommand{\bLambda}{\boldsymbol{\Lambda}}
\newcommand{\bXi}{\boldsymbol{\Xi}}
\newcommand{\bSigma}{\boldsymbol{\Sigma}}
\newcommand{\bUpsilon}{\boldsymbol{\Upsilon}}
\newcommand{\bPhi}{\boldsymbol{\Phi}}
\newcommand{\bPsi}{\boldsymbol{\Psi}}
\newcommand{\bOmega}{\boldsymbol{\Omega}}

\def\theequation{\arabic{section}.\arabic{equation}}
\setcounter{page}{1}
Study of instantaneous dependence among several variable is important in many of the high-dimensional sciences. Multivariate GARCH models are as a standard approach for modelling time-varying covariance matrix such phenomena. Cholesky GARCH is one of these approaches where the time-varying covariance matrix can be written parsimoniously, containing variance components through a diagonal matrix and dependency components through a unit lower triangular matrix with regression coefficients as entries. In this paper, we proposed a stochastic structure for dependency components in Cholesky GARCH model by considering linear regression model as a state-space model and using kalman filtering for estimating regression coefficients. We find that the MSE of stochastic Cholesky GARCH model is smaller than the MSE of other models also show that the stochastic Cholesky GARCH has better performance in compare to another models based on MAE and  MSE criterions for the real data.

\section{Introduction}
Accurate estimation of the covariance matrices, as well as the time-varying covariance matrices, is important in many of the high-dimensional sciences to study relationship between variables and gain time-varying correlation matrices. In these sciences, require the modeling and forecasting of time-varying covariance matrices $\Sigma_t$ based on the independently $N(0,\Sigma_t)$-distributed data
$Y_t$ , $t = 1, 2, . . . , n$, where
$Y_t$ is the $p$-dimensional variable at time $t$ of a multivariate time
series.
There are two important challenges in modelling covariance matrices. The first problem is the curse of dimensionality, in a $p$-dimensional time series, the time-varying covariance matrix $y_t$ consists of $\frac{p(p+1)}{2}$ different time-varying elements. The second problem is maintaining the positive-definite constraint. The  time-varying covariance matrix $\Sigma_t$ must be positive definiteness for all t. Special attention is needed to  maintain the constraint when  N is large. 
Dynamic dependence $\Sigma_t$ is the subject of multivariate volatility modeling. A variety of multivariate extensions of the univariate generalized autoregressive conditional heteroscedastic (GARCH) models (Bollerslev, 1986 \cite{bollerslev:1986}) has been developed in the finance literature. VEC and BEKK models are the first class of multivariate GARCH models which arised as a direct generalization of the univariate GARCH models. At these models, number of parameters increase rapidly with dimensionality and parameter estimation is heavy computationally. In some other existing methods impose strong assumptions on the conditional correlation matrices. For instance, Bollerslev (1990) \cite{bollerslev:1990} assumes constant conditional correlation for the GARCH model (CCC-GARCH), which may not be satisfied in real data. Engle (2002) \cite{engle:2002} and Tse and Tsui (2002) \cite{tse etal:2002} assume the dynamic conditional correlation for the GARCH model (DCC-GARCH), which is computationally expensive in high-dimensional cases. To overcome the difficulty of curse of dimensionality, several methods have been proposed in the literature for multivariate volatility modeling. In
particular, the idea of orthogonal transformation has attracted much attention. For instance, Alexander (2001) proposed Orthogonal-GARCH (O-GARCH) where she has been used principal component analysis (PCA) in multivariate
volatility modelling. Also Van derWeide (2002) extracts the concept of ICA to present
a class of generalized orthogonal GARCH (Go-GARCH) models for volatility
modelling. Matteson and Tsay (2011) propose a dynamic orthogonal component (DOC)
method to multivariate volatility modelling where the components are uncorrelated.\\
The modified Cholesky decomposition (MCD) of a covariance matrix is another the orthogonal transformation approach. The use of MCD to model multivariate volatility, guarantees positive definiteness volatility matrices. In addition, it provides a parameterization of the covariance matrix with unconstrained and statistically interpretable parameters (pourahmadi 1999 \cite{pourahmadi:1999}). Using the MCD, the covariance matrix can be decomposed parsimoniously, with variance components through a diagonal matrix and dependency components through a unit lower triangular matrix with $\frac{p(p-1)}{2}$  elements.\\
 Recently, Dellaportas et al. (2012) \cite{dellaportas et al:2012}, Pedeli et al. (2015) \cite{pedeli et al:2015} considered multivariate volatility matrix and proposed  Cholesky GARCH (C-GARCH) and Cholesky Log-GARCh (CLGARCH). They assumed elements of the lower triangular matrix of Cholesky decomposition to be constant over time. Therefore they considered a constant dependence structure and time-varying nature of the diagonal matrices.\\ 
 In this paper, we propose a time-varying dependence structure for MCD. Thus for the estimation of corresponding linear regression, we consider some state space model and estimate its parameters by using kalman filter. This method updates information on regression coefficients by implying informations as arrive. In result of, the regression coefficients are updated through the time. So, we make a time-varying covariance matrix base on MCD idea by a time-varying dependence structure and a  time-varying variance structure as GARCH models.\\
 This article is organized as follows:
we review the MCD approach to model a covariance matrix in section 2. Section 3
explains time-varying Cholesky factor based on state-space models. In section 4 is about Flexible cholesky GARCH model and addressing estimation and order issue.  Section 5 explains simulation study and demonstrates its result. Application are illustrated in section 6.

\section{Modified Cholesky Decomposition}
This section reviews the role of the modified cholesky decomposition in reparametrizing a covariance matrix in term of unconstrained regression coefficients which guarantees the positive definiteness of the estimated covariance matrix.\\
let $Y=(Y_1,..,Y_p)^{\prime}$ is a $p$-dimensional vector of mean zero random variables with the covariance matrix $\Sigma$. This decomposition arise from regressing each variable $Y_j$ on its predecessors $Y_{j-1},...,Y_1$, for $j=2,...,p$ is defined as

\begin{equation}
         Y_j=\sum_{k=1}^{j-1}\phi_{jk}Y_k+\epsilon_j  \label{eq.MCD1}
 \end{equation}
where $\phi_{j1},...\phi_{j(j-1)}$ are regression coefficients and $\epsilon_j$ denotes the linear least-squares prediction error with variance $\sigma_j=var(\epsilon_j)$ and $\epsilon_1=Y_1$. Let $\epsilon=Y-\hat{Y}=(\epsilon_1,...,\epsilon_p)$ be the vector of successive uncorrelated prediction errors with diagonal covariance matrix 

\begin{equation*}
        cov(\epsilon)=diag(\sigma^2_1,...\sigma^2_p)=D
 \end{equation*}
The \ref{eq.MCD1} can be written in matrix form as

\begin{equation}\label{form matrix}
       TY=\epsilon
 \end{equation}
where $T$ is a unite lower triangular matrix with the negative of the regression coefficients $\phi_{jk}$ as it's entries:
\begin{equation}
T = \left( {\begin{array}{*{20}{c}} \label{T}
1&0&0& \cdots &0\\
{{-\phi _{21}}}&1&0& \cdots &0\\
{{-\phi _{31}}}&{{-\phi _{32}}}&1& \cdots &0\\
 \vdots & \vdots & \ddots & \vdots & \vdots \\
{{-\phi _{p1}}}&{{-\phi _{p2}}}& \cdots &{{-\phi _{p,p - 1}}}&1
\end{array}} \right)
\end{equation}

Then, the MCD of $\Sigma$ follow from
\begin{equation}
D=cov(\epsilon)=cov(TY)=Tcov(Y)T^{\prime}=T\Sigma T^{\prime}
\end{equation}
Or equivalently 
\begin{equation}
\Sigma=T^{-1}D T{^{\prime}}^{-1}
\end{equation}
Thus, the MCD of a covariance matrix provides a parameterisation of the covariance matrix with unconstrained parameters, and converts the difficult task of modelling a covariance matrix to that of modelling te sequence regression in \ref{eq.MCD1}.
\setcounter{equation}{0}

\section{Time-varying Cholesky factor based on state-space models}
Let for a given time point $t$, $\mathbf{Y}_t= ({Y}_{1t},{Y}_{2t},\dots,{Y}_{pt})$ be a $p$-dimensional vector of mean zero Gaussian random variables and covariance matrix $\mathbf{\Sigma}_t$. The linear regression $Y_{jt}$ on its predecessors, $Y_{1t},Y_{2t},\dots,Y_{(j-1)t}$ is defined so
\begin{equation}\label{reg}
Y_{jt}=\sum_{k=1}^{j-1}\phi_{jk}Y_{kt}+\epsilon_{jt} \quad t=1,2,\ldots,n,\quad  j=2,...,p
\end{equation}
where $\phi_{kt}$'s are regression coefficients and
 $\epsilon_{jt}$ are independent Gaussian with mean zero and variance $\sigma^2_{jt}$.\\
Now, we assume that the parameters in \ref{reg} are time dependent, where the parameters are allowed to evolve through time according to multivariate Gaussian random walk. So the equation \ref{reg} can be written as 
\begin{equation}\label{time-reg}
Y_{jt}=\mathbf{X^{\prime}}_{jt}\mathbf{\phi}_{jt}+\epsilon_{jt} \quad t=1,2,\ldots,n,\quad  j=2,...,p
\end{equation}
where $\mathbf{\phi}_{jt}=(\phi_{j1t},\ldots,\phi_{j(j-1)t})$ and $\mathbf{X}_{jt}=(Y_{1t},Y_{2t},\ldots,Y_{(j-1)t})^{\prime}$ is 
the \textit{t}-th
observed vector of ${Y}_{2t},{Y}_{3t},\dots,{Y}_{(j-1)t}$.
The state-space representation \ref{time-reg} is written as
\begin{equation} \label{state}
Y_{jt}=\mathbf{X^{\prime}}_{jt}\mathbf{\phi}_{jt}+\epsilon_{jt}
\end{equation}
$$\mathbf{\phi}_{jt}=\mathbf{A}_j\mathbf{\phi}_{j(t-1)}+\mathbf{a}_{jt}$$
where $\mathbf{A}_j$ is a $(j-1)\times (j-1)$ identity matrix and $\mathbf{a}_{jt}=(a_{j1t},...,a_{j(j-1)t})$ are white noises at time $t$, with covariance matrix $\mathbf{Q}$. The first equation is an observation equation and the second equation the state equation, which describes the evolution of the state vector.\\
  To solve the recursive equation in \ref{state} by Bayesian method, it is necessary to have inferences about the initial values of the state vector $\mathbf{\phi}_{jt}$ and its distribution. The probability distribution of state vector at time zero is denoted by $p(\mathbf{\phi}_{j0})$. The state equation \ref{state} with such initial distributions determines the distribution of the state vectors $p(\mathbf{\phi}_{jt}|\mathbf{\phi}_{j(t-1)})$, $t=1,\ldots,n$ and $j=2,\ldots,p-1$. We refer to these distributions as prior distributions of state vector at time $t$. To specify the likelihood model of measurements as $p(y_{jt}|\mathbf{\phi}_{jt})$, after observing a new measurement $y_{jt}$, the prior distributions to are updated, such update is denoted by $\mathbf{\phi}_{jt}$ expressed by $p(\mathbf{\phi}_{jt}|y_{jt})$, called posterior distribution. By the assumption that the errors $\epsilon_{jt}$ and $\mathbf{a}_{jt}$, are normality distributed we conclude that as 
   \begin{equation}
   p(Y_{jt}|\mathbf{\phi}_{jt})=N(Y_{jt}|\mathbf{X^{\prime}}_{jt}\mathbf{\phi}_{jt},\sigma^2_{jt})
   \end{equation}
   $$p(\mathbf{\phi}_{jt}|\mathbf{\phi}_{j(t-1)})=N(\mathbf{\phi}_{jt}|\mathbf{\phi}_{j(t-1)},\mathbf{Q})$$
   where $\mathbf{Q}$ is covariance matrix of the multivariate random walk noises. Now, given the distribution
   $$p(\mathbf{\phi}_{j(t-1)}|Y_{j(1:t-1)})=N(\mathbf{\phi}_{j(t-1)}|\mathbf{\phi}_{j(t-1)},\mathbf{P}_{j(t-1)})$$ 
   the joint distribution of $\mathbf{\phi}_{jt}$ and $\mathbf{\phi}_{j(t-1)}$ when we have the information of $Y_{j(1:t-1)}$ can be expressed as
   $$p(\mathbf{\phi}_{jt},\mathbf{\phi}_{j(t-1)}|Y_{j(1:t-1)})=p(\mathbf{\phi}_{jt}|\mathbf{\phi}_{j(t-1)})p(\mathbf{\phi}_{j(t-1)}|Y_{j(1:t-1)})$$
   So,
   $$p(\mathbf{\phi}_{jt}|Y_{j(1:t-1)})=\int p(\mathbf{\phi}_{jt}|\mathbf{\phi}_{j(t-1)})p(\mathbf{\phi}_{j(t-1)}|Y_{j(1:t-1)}) \,d\mathbf{\phi}_{j(t-1)}$$
  As $p(\mathbf{\phi}_{jt}|\mathbf{\phi}_{j(t-1)})$ and $p(\mathbf{\phi}_{j(t-1)}|Y_{j(1:t-1)})$ are Gaussian, the result of marginalization is Gaussian 
  $$p(\mathbf{\phi}_{jt}|y_{j(1:t-1)})=N(\mathbf{\phi}_{jt}|\mathbf{\phi}_{jt}^-,\mathbf{P}_{j(t)}^-)$$
  where 
  \begin{equation}\label{prediction}
  \mathbf{\phi}_{jt}^-=\mathbf{A} \mathbf{\phi}_{j(t-1)}
  \end{equation}
  $$\mathbf{P}_{j(t)}^-=\mathbf{A}\mathbf{P}_{j(t-1)}\mathbf{A}^{\prime}+\mathbf{Q}$$
  By using this as the prior distribution for the measurement $p(y_{jt}|\mathbf{\phi}_{jt})$ we gain posterior distribution as
  $$p(\mathbf{\phi}_{jt}|y_{j(1:t)})=N(\mathbf{\phi}_{jt}|\mathbf{\phi}_{jt},\mathbf{P}_{j(t)})$$
  where the Gaussian parameters can be obtained by completing the quadratic form in the exponent, which gives
  \begin{equation}\label{phi}
  \mathbf{\phi}_{jt}=\left[ (\mathbf{P}_{j(t)}^-)^{-1}+\frac{1}{\sigma_j^2}\mathbf{x^{\prime}}_{jt}\mathbf{x}_{jt}\right]^{-1}\left[ \frac{1}{\sigma_j^2}\mathbf{x^{\prime}}_{jt}y_{jt}+(\mathbf{P}_{j(t)}^-)^{-1}\mathbf{\phi}_{jt}^-\right]
  \end{equation}
  $$  \mathbf{P}_{j(t)}=\left[(\mathbf{P}_{j(t)}^-)^{-1}+\frac{1}{\sigma_j^2}\mathbf{x^{\prime}}_{jt}\mathbf{x}_{jt}\right]^{-1 }$$
  This recursive computational algorithm for estimating parameters is a special case of the Kalman filter algorithm.\\
  Now, we can write matrix form in \ref{form matrix} as
  \begin{equation}\label{form matrix t}
  T_tY_t=\epsilon_t
  \end{equation}
    where $T_t$ is a time-varying matrix $T_t$ with
  the negative of the regression coefficients $\phi_{jkt}$ that obtained in \ref{phi} as it's entries
 \begin{equation}
T_t = \left( {\begin{array}{*{20}{c}} \label{timeT}
1&0&0& \cdots &0\\
{{-\phi _{21t}}}&1&0& \cdots &0\\
{{-\phi _{31t}}}&{{-\phi _{32t}}}&1& \cdots &0\\
 \vdots & \vdots & \ddots & \vdots & \vdots \\
{{-\phi _{p1t}}}&{{-\phi _{p2t}}}& \cdots &{{-\phi _{p,(p - 1)t}}}&1
\end{array}} \right).
\end{equation} 
   
So, for each time $t$, we will have MCD covariance matrix as $\Sigma_t=T_t^{-1}D_t T_t{^{\prime}}^{-1}$.

\section{Stochastic Cholesky GARCH} 
From the matrix form of \ref{form matrix t} we have 
\begin{equation}
T_tY_t=\epsilon_t\sim N(0,D_t)
\end{equation}
where $D_t=diag(\sigma_{1t}^2,...,\sigma_{pt}^2)$. In this paper, for each time -varying innovation variance $\sigma_{jt}^2$, $j=1,...,p$, we model $\sigma_{jt}^2$ by a GARCH(p,q) model defined recursively in time as,
\begin{equation}\label{GA}
 \sigma^2_{jt}=w+\sum_{i=1}^p\alpha_{i}\epsilon_{t-i}^2+\sum_{l=1}^q\beta_{l}\sigma^2_{j(t-l)} 
\end{equation}
In the following, we pay to estimation parameters of time-varying covariance matrix $\Sigma_t$ and issue of ordering the stocks a portfolio which is a important object in our method.
\subsection{Estimation}
Estimation of time-varying covariance matrix $\Sigma_t$ on based Cholesky GARCH, are performed in two steps. At first, we should estimate time-varying cholesky factor $T_t$ as explained in section $2$ and then, using $\phi_{jkt}$'s, vector prediction errors (innovation) be obtain as 
\begin{equation}
\epsilon_{jt}=\begin{cases}
Y_{1t}, & j=1\\
Y_{jt}-\mathbf{X^{\prime}}_{jt}\mathbf{\phi}_{jt}, & j=2,...p
\end{cases}
\end{equation}
\\
In the second step, to estimate time-varying diagonal matrix $D_t$, we need to estimate parameters in GARCH model \ref{GA}. For that, assuming normality for the returns, log-likelihood function, with ignoring an irrelevant constant is given by

     \begin{align*}
           L(\mathbf{\theta}) & =\sum_{t=1}^n (log\vert\Sigma_t\vert+Y_t^{\prime}\Sigma_t^{-1}Y_t) \\
           & =\sum_{t=1}^n(\sum_{j=1}^p log\sigma_{jt}^2+Y_t^{\prime}T_t^{\prime}D_t^{-1}T_tY_t)\\
            & =\sum_{t=1}^n\sum_{j=1}^p(log \sigma_{jt}^2+\frac{\epsilon_{jt}^2}{\sigma_{jt}^2})\\
            &=\sum_{j=1}^p\left\lbrace \sum_{t=1}^n\left( log \sigma_{jt}^2+\frac{\epsilon_{jt}^2}{\sigma_{jt}^2}\right)\right\rbrace
            \end{align*}
 Where $\sigma_{jt}$ is a GARCH model as defined at \ref{GA} and $\mathbf{\theta}$ is vector of parameters in a cholesky GARCH. If we focus on the case of a model of order (1,1), vector of parameters $\mathbf{\theta}$ is as $\mathbf{\theta}=(w,\alpha_{1},\beta_{1})$.
 To obtain vector parameters $\theta$, the log-likelihood function can be maximized by standard numerical procedures. So now, we can estimate  time-varying diagonal matrix $\hat{D_t}$ for each time $t$.\\
 Finally, be estimated time-varying matrix $\Sigma_t$ for each time $t$ as following
 \begin{equation}
  \hat{\Sigma_t}=\hat{T_t}^{-1}\hat{D_t}\hat {T_t}{^{\prime}}^{-1}
 \end{equation}
 Ordering the stocks in a portfolio is a drowback in GARCH models based on idea MCD which is discussed next.

\subsection{Ordering}
Ordering variables is a great challenge in statistical desicion problem and a long-standing problem on statistics. For a portfolio of $p$ stocks, there are $p!$ choices.\\
Basford and Tukey (1999) propose a method for ordering variables based on a scatterplot matrix nice in the sense that one brings the more correlated variables closer to the main diagonal, challed "greedy close algorithm". On the basis "greedy close algorithm", Bickel and levina (2008) and Bickel and Gel (2011) provided idea more powerful perceptual idea of "banded sample covariance matrix estimation". In the context order variables in regresion, Didge and Rousson (2001) present a ordering variables in setting regression, just for two variables, based on asymmetric properties correlation coefficient. Dellaportas and Pourahmadi
(2004; 2012) suggested using Akaike information criterion (AIC) or Bayesian information criterion (BIC) for ordering variables in MCD. Also, Pedeli and et al. (2015) used the same criterion in their paper. \\
Kang and et al. (2016) suggest a novel order-averaged MCD-based approach for estimating
covariance matrices using a random sample from the population of $p!$ permutations of
the variables. 
In this paper, we use of BIC criterion for ordering variables in MCD.

\section{Simulation}
In this section we conduct two simulation experiments to investigate the consistency estimators and  evaluate the performance of the different methods for estimating time-varying
correlation matrix.
\subsection*{Simulation 1}
In this simulation we investigate the consistency parameters $\mathbf{\phi}_{jt}$ and parameters in GARCH model.
Four  sample lengths n=100, 500, 1000 and 2000 observations have been used in experiment, and there are 1000 replications for each sample size.  For each time point, we have one estimate for parameter $\phi_{jt}$, so for shorthand, we just present  ten of $\phi_{jt}$ for each sample size. Table \ref{tab1-s} shows  bias values and estimate of parameters $\phi_{jt}$ for ten time point. It can be seen that  Bias is generally small and decrease as the sample size increases. 

  \begin{table}[h]
\caption{} \label{tab1-s}
\centering
\small
 \begin{tabular}{ccccccccc}
 \hline 
\hline 
 & $n=100$& &$n=500$ && $n=1000$ 
\\\hline 
parameter&  real value & estimate&  real value & estimate& real value & estimate&
\\\hline
$\phi_{1(n-9)}$ & 0.4145197946    &0.4680577377 &-0.2661494  &  -0.3024227&0.6787251&          0.7170813 \\
$\phi_{1(n-8)}$ & 0.4391433832&0.4613561505& -0.2338428  & -0.2532158  & 0.6564266 &  0.6411753        \\
$\phi_{1(n-7)}$ &  0.3925500472 &0.3751479629  &-0.3241438 & -0.3906113 &0.6698831& 0.6625149        \\
$\phi_{1(n-6)}$ &  0.3344852256&0.3880605372 &   -0.2832782 &-0.3090194&0.6069042 & 0.6246090      \\
$\phi_{1(n-5)}$ &  0.3608330633 & 0.3636095009& -0.3273064  &-0.3612746 &  0.5984047 &0.6246090           \\
$\phi_{1(n-4)}$ &  0.3301939939  & 0.3665939630 & -0.3351310   & -0.3682334 & 0.6422929 & 0.6695019     \\
$\phi_{1(n-3)}$&  0.3050611094  & 0.2891016385  &-0.3378436  &-0.2948399 &0.6380926 & 0.7172956         \\
$\phi_{1(n-2)}$&   0.2872172739 & 0.2758778669&-0.3789086 &-0.4160425  & 0.6441557& 0.5484562       \\
$\phi_{1(n-1)}$&   0.3105050763& 0.2995623391  &-0.4440797  &-0.3903170&0.6550884&0.7163191      \\
$\phi_{1n}$     & 0.3135250794 & 0.2696930218   & -0.4766288 &-0.4676340&0.6262435&0.5614452      \\
\\\hline
Bias &-0.009317485&&0.006857815&&-0.0006963098
\\\hline
\end{tabular}
\end{table}

\subsection*{Simulation 2}
\begin{figure}[ht]
\begin{center}
\includegraphics[width=12cm,height=13cm]{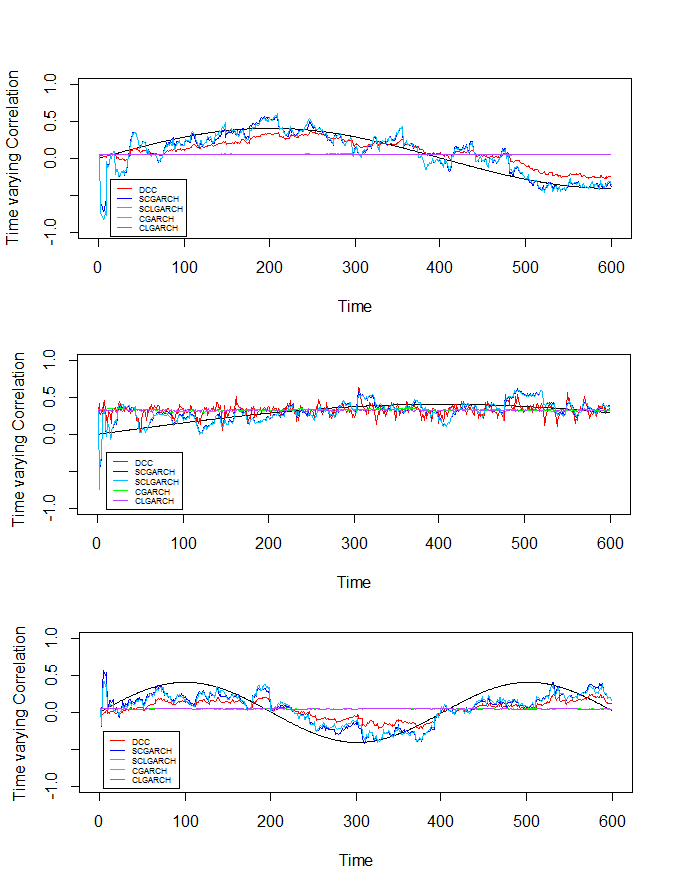}
\caption{\footnotesize{\textbf Result fitting DCC, CGARCH, CLGARCH, SCGARCH and SCLGARCH models, top chart shows time varying correlation related to $\sigma_{21t}$, middle chart shows time varying correlation related to $\sigma_{31t}$ and lower chart shows time varying correlation related to $\sigma_{32t}$.}}\label{fig1}
\end{center}
\end{figure}

 In this simulation, we consider $p=3$. So, we generate $\mathbf{y_t}$ from bivariate normal distribution with mean zero and covariance matrix as
\begin{eqnarray}
	\label{eq14}
\mathbf{\Sigma_t}=\left( {\begin{array}{*{20}{c}} \label{T}
2&{{\sigma _{21t}}}&{{\sigma _{31t}}}\\
{{\sigma _{21t}}}&3&{{\sigma _{32t}}}\\
{{\sigma _{31t}}}&{{\sigma _{32t}}}&4\\
\end{array}} \right)
\end{eqnarray}
where the covariance term $\sigma_{21t}$,  $\sigma_{31t}$ and  $\sigma_{32t}$ varies across time for $t=1,...,n$. This
allowed us to control the dynamic relationship between the two time courses $y_{1t}$, $y_{2t}$ and $y_{3t}$
throughout the time series. we suppose that connectivity of between
brain regions is in the form of sine. The value of $\sigma_{21t}$, $\sigma_{31t}$ and $\sigma_{32t}$ were set as $\sigma_{21t}=sin(t/\delta), \delta=\frac{1024}{2^3}$, $\sigma_{31t}=sin(t/\delta), \delta=\frac{1024}{2^2}$ and $\sigma_{32t}=sin(t/\delta), \delta=\frac{1024}{2^4}$, respectively. \\
 Simulation study was repeated 1000 times, and the MSE of five models represented in Table \ref{tab1}.  As it is seen, in Table \ref{tab1}, value of MSE SCGARCH model is smaller MSE other models. This shows that the SCGARCH model has a better performance than other models in this simulation.  Figure\ref{fig1} shows  the results of fitting DCC, CGARCH, CLGARCH, SCGARCH and SCLGARCH models in one of the
iterations. \\
 \begin{table}[h]
\caption{The results of the MSE simulation for DCC, CGARCH, CLGARCH, SCGARCH and SCLGARCH models} \label{tab1}
\centering
\small
 \begin{tabular}{cccccccc}
\hline 
Model & MSE
\\\hline
DCC-GARCH & 0.01260446  \\
CGARCH& 0.006155253   \\
CLGARCH &  0.006646985 \\ 
SCGARCH &0.08321035  \\
SCLGARCH &0.08358656\\
\hline
\end{tabular}
\end{table}

\section{Application}
In order to examine and evaluate the ability of the proposed model to estimate the instantaneous relationship between variables, we survey two sets of data from two different scientific fields about neuroscience and finance.\\
The first data set related to functional magnetic resonance 
imaging (fMRI) from brain regions. One of the main studies areas in fMRI is describing the functional connectivity 
(FC) between brain regions. This study is followed by describing the dependency structure among of blood-oxygen-level-dependent (BOLD) response time series from various brain regions.\\ Previously, most fMRI studies have implicitly
assumed these time series to be
stationary over time. 
Recently, scientists
have observed indications of non-stationarity in these
time series.This led to an increased interest in attempting
to quantify the dynamic changes in functional connectivity (FC) during the course of an fMRI experiment. Here, we use multivariate time series for obtaining dynamic FC and compare the results together (see Lindquist et. al. (2014)).\\
Second data set related to a stock portfolio in finance. Many tasks of financial management, such as stock portfolio selection, option pricing and risk assessment, require the modelling and forecasting of time-varying correlation matrices.\\ This stock portfolio consists of several stocks that we would like to get their instant connection by using multivariate time series and compare the results together.
To measure the accuracy of a covariance matrix estimate $\hat{\mathbf{\Sigma}}_t$, we consider as a
benchmark the moving blocks approach (i.e., $\mathbf{\Sigma_t}=\tilde{\mathbf{\Sigma}}_t$) and calculate the mean absolute error and mean squared error given by:
$$MAE_t=\frac{1}{p^2}\sum_{i=1}^{p}\sum_{j=1}^{p}|\hat{\sigma}_{ijt}^2-\sigma_{ijt}^2|$$
$$MSE_t=\frac{1}{p^2}\sum_{i=1}^{p}\sum_{j=1}^{p}(\hat{\sigma}_{ijt}^2-\sigma_{ijt}^2)^2$$
respectively. For the comparison of different methods we rely on the averages of the
above measures over $t$.\\
A notable challenge in computing any such measure of accuracy is that the true covariance
matrix $\mathbf{\Sigma_t}$ is unknown. We resolve this challenge by employing a moving block
technique to get a reliable proxy for it Lopes et al. (2012).
Selection of the block size $q$ was based on MAE and MSE criteria with the moving
blocks estimator serving as $\hat{\mathbf{\Sigma}}_t$ and $\mathbf{\Sigma_t}$ being the observed covariance matrix. More specifically, the average loss functions $MAE=\sum_{t=1}^{n}MAE_{t}/n$ and $MSE=\sum_{t=1}^{n}MSE_{t}/n$ were calculated
for several values of $q$. Let $q_1,...,q_k,...,q_p$ be the ordered series of $q$-values.
The value $q_k$ that stabilizes the absolute difference 
$|\hat{\bigtriangleup}^{(q_k)}-\hat{\bigtriangleup}^{(q_{k-1})}|$ for two loss functions under consideration, where $\hat{\bigtriangleup}$ denote MAE and MSE, was selected as the optimal $q$. According to this procedure, we chose $q=35$, $q=65$ in first data set and second data set, respectively.

\subsection{fMRI data experiment based on default mode network (DMN) regions }
\begin{figure}[ht]
\begin{center}
\includegraphics[width=13cm,height=5cm]{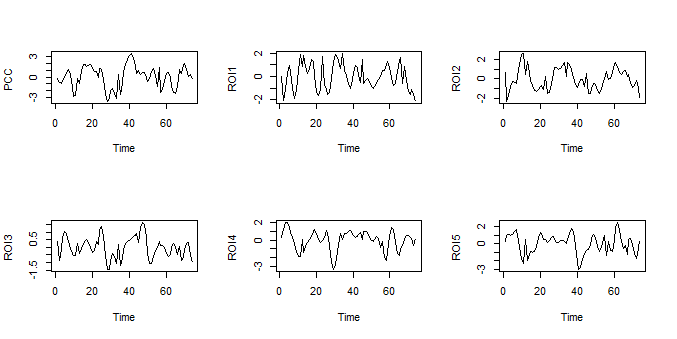}
\caption{\footnotesize{Time series plots of brain regions from default mode network (DMN).  }}\label{fig fmri}
\end{center}
\end{figure}

\begin{figure}[ht]
\begin{center}
\includegraphics[width=12cm,height=8cm]{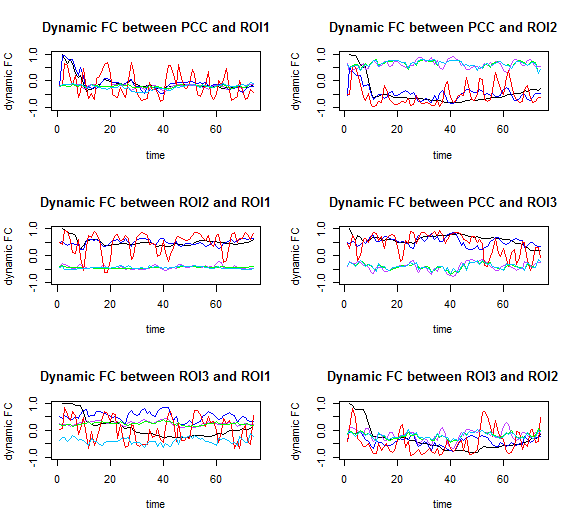}
\includegraphics[width=12cm,height=8cm]{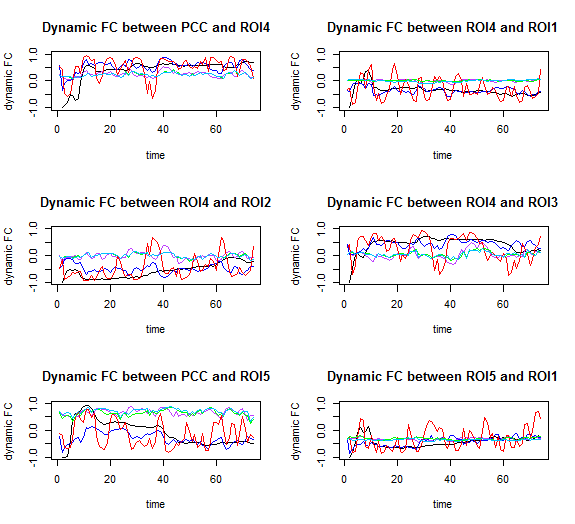}
\caption{\footnotesize{Estimate dynamic FC based on DCC GARCH (red curve), CGARCH (green curve), CLGARCH (dark orchid curve), SCGARCH (blue curve), SCLGARCH  (deep sky blue curve) models and true dynamic correlation that obtained based on moving block approach  (black curve) among regions.}}\label{resultfmri1}
\end{center}
\end{figure}

\begin{figure}[ht]
\begin{center}
\includegraphics[width=12cm,height=8cm]{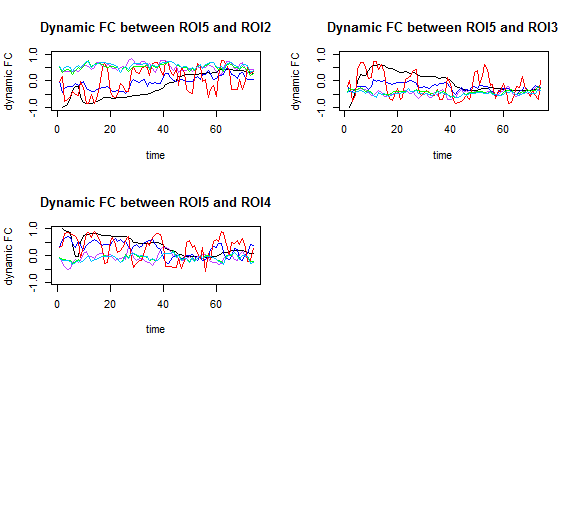}
\caption{\footnotesize{Estimate dynamic FC based on DCC GARCH (red curve), CGARCH (green curve), CLGARCH (dark orchid curve), SCGARCH (blue curve), SCLGARCH  (deep sky blue curve) models and true dynamic correlation that obtained based on moving block approach  (black curve) among regions.}}\label{resultfmri3}
\end{center}
\end{figure}

In this study, we used fMRI data on brain images of a global competition called ADHD-200, which was held in 2011, whose pre-processed collection was compiled in a resource called ADHD-200 Preprocessed and placed at $www.nitrc.org$ address.
These brain images are about resting-state, from these images, we selected an image randomly. As mentioned, these data are pre-processed, so we don't need to do pre-processing stage of data. For extracting brain regions time
series, used $wfu-pickatlas$ toolbox in MATLAB. \\
Brain regions of interest (ROIs) were selected based on Chang et.al. (2010). In the default mode network (DMN), the posterior
cingulate cortex (PCC) was selected as the principal region and
five other regions of interest (ROIs), demonstrated by  Chang et.al. (2010)
as brain regions with the most enormous variation in
dynamic correlation with PCC, were also chosen. These
brain regions are as, $ROI1$ Right inferior parietal cortex, $ROI2$ Right inferior frontal operculum, $ROI3$ Right inferior temporal cortex, $ROI4$ Right inferior orbitofrontal cortex, $ROI5$ Anterior cingulate cortex (ACC).\\
Time series plots of brain regions from default mode network is showed in figure \ref{fig fmri}. Table \ref{tab2} summarizes the results of comparing the performance of three methods for experimental data. This table shows that our proposed model outperforms other models in terms of two measures of accuracy. Figure \ref{resultfmri1} and Figure \ref{resultfmri3} show the result of the estimate dynamic FC based on our proposed model, DCC-GARCH model, Cholesky GARCH and true dynamic correlation that obtained based on moving block approach among regions.

  \begin{table}[h]
\caption{The results of the MAE and MSE for DCC, CGARCH, CLGARCH, SCGARCH and SCLGARCH models in fMRI data. } \label{tab2}
\centering
\small
 \begin{tabular}{cccc}
\hline 
Method & MAE  & MSE
\\\hline
DCC-GARCH &  0.2880811  &0.1673271 \\
CGARCH&     0.4400117          &   0.3752082         \\
CLGARCH&  0.4384204          &  0.3704371                \\
SCGARCH&     0.1585642       &     0.05777206           \\
SCLGARCH&  0.4477702       &     0.3827415            \\
\hline
\end{tabular}
\end{table}

\subsection{Monthly stock return of 12 US bluechips}
\begin{figure}[ht]
\begin{center}
\includegraphics[width=14cm,height=9cm]{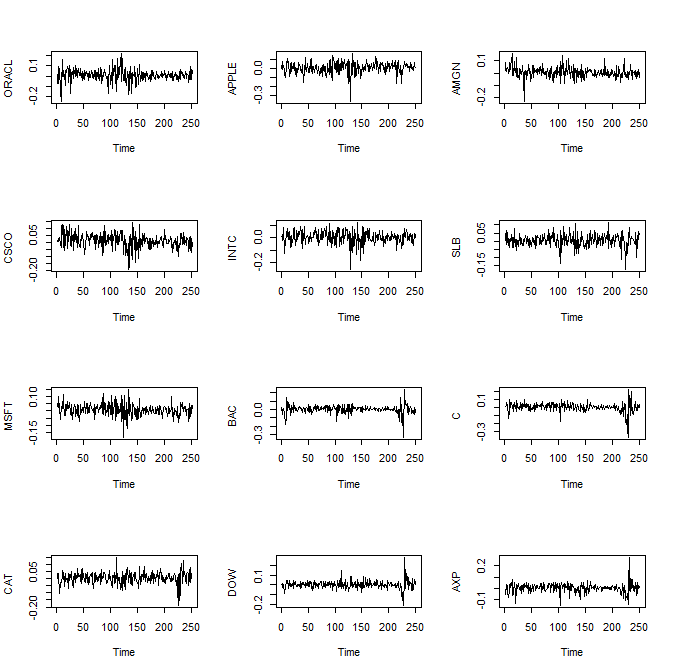}
\caption{\footnotesize{Time series plots of monthly returns of 12 US bluechips for the period January 1990 to December 2010.}}\label{fig time12 us}
\end{center}
\end{figure}

To compare the performance of the proposed model with other models, we consider compound monthly log-returns from a panel of 12 US bluechips .These 12 US bluechips includes companies Oracle Corp. (ORCL), Apple Inc.
(AAPL), Amgen Inc. (AMGN), Cisco Systems Inc. (CSCO), Intel Corp. (INTC), Schlumberger Ltd. (SLB), Microsoft Corp. (MSFT), Bank of America Corp. (BAC), Citigroup Inc. (C), Caterpillar Inc. (CAT), Dow Chemical Co. (DOW) and American Express Co. (AXP). The n = 251 stock returns were collected
from January 1990 to December 2010. Time series plots of monthly log-returns this portfolio is showed in figure \ref{fig time12 us}. Table \ref{tab3} summarizes the results of comparing the performance of three methods for 12 US bluechips data from January 1990 to December 2010.

  \begin{table}[h]
\caption{The results of the MAE and MSE for DCC, CGARCH, CLGARCH, SCGARCH and SCLGARCH models in stock returns data} \label{tab3}
\centering
\small
 \begin{tabular}{cccc}
\hline 
Method & MAE  & MSE
\\\hline
DCC-GARCH &  0.1081201 &0.02267042 \\
CGARCH&   0.2087609         &   0.07335637        \\
CLGARCH&    0.204851      & 0.07028564             \\
SCGARCH& 0.09565078         &   0.01899578       \\
SCLGARCH&   0.1059657  &     0.02249688        \\
\hline
\end{tabular}
\end{table}




\end{document}